%% file: main.tex
\documentclass[10pt]{amsart}
\usepackage{import}
\usepackage{./preamble/Pakete}
\usepackage{./preamble/Befehle}

\begin{document}
\import{Inhalt/}{title}
\section{Introduction}
\import{Inhalt/}{intro.tex}

\section{Motivic sheaves d’apr\`es Cisinski--D\'eglise} 
\import{Inhalt/}{motives.tex}
\section{Motivic sheaves on affinely stratified varieties}
\import{Inhalt/}{affinelystratified.tex}
\section{Flag varieties and Koszul duality}
\import{Inhalt/}{flagvariety.tex}
\appendix
\section{Weight structures and t-structures}
\import{Inhalt/}{tandweight.tex}
\bibliographystyle{amsalpha} 
\bibliography{main}

\end{document}

%% file: Inhalt/title.tex

\title[$K$-Motives and Koszul duality]{$K$-Motives and Koszul duality}

\author{Jens Niklas Eberhardt}
\address{Max-Planck-Institut fuer Mathematik\\ Vivatsgasse 7, 53111 Bonn, Germany
}
\email{mail@jenseberhardt.com}  
\begin{abstract}
We construct an \emph{ungraded version} of Beilinson--Ginzburg--Soergel's Koszul duality for Langlands dual flag varieties, inspired by Beilinson's construction of rational motivic cohomology in terms of $K$-theory. 

For this, we introduce and study categories $\KM_\Ss(X)$ of $\Ss$-constructible $K$-motivic sheaves on varieties $X$ with an affine stratification. We show that there is a natural and geometric functor, called \emph{Beilinson realisation}, from $\Ss$-constructible mixed sheaves $\Dmix_\Ss(X)$ to $\KM_\Ss(X).$ 

We then show that Koszul duality intertwines the Betti realisation and Beilinson realisation functors and descends to an equivalence of constructible sheaves and constructible $K$-motivic sheaves on Langlands dual flag varieties. \end{abstract}
\maketitle
\setcounter{tocdepth}{1} 
\tableofcontents

%% file: Inhalt/intro.tex
\todo{A1 (additional change 1)}
Let $G\supset B$ be a split reductive group with a Borel subgroup and $X=G/B$ be the flag variety. Denote by $G^\vee\supset B^\vee$ the Langlands dual group with a Borel subgroup and by $X^\vee=G^\vee/B^\vee$ the Langlands dual flag variety. \todo{S1 (specific comment 1)}

In this article we prove an \emph{ungraded} and \emph{$K$-theoretic} version of Beilinson--Ginzburg--Soergel's Koszul duality.
\begin{theorem*} There is a commutative diagram of functors where the horizontal arrows are equivalences
\begin{center}
\begin{tikzcd}
\Dmix_{(B)}(X) \arrow{r}{\gKos} \arrow{d}{v}
& \Dmix_{(B^\vee)}(X^\vee) \arrow{d}{\iota} \\
\D_{(B)}(X)  \arrow{r}{\Kos}
& \KM_{(B^\vee)}(X^\vee).
\end{tikzcd}\end{center}
\end{theorem*}
Let us explain the ingredients of this diagram.  In \cite{BG}, \cite{Soe90} and \cite{BGS}, Beilinson--Ginzburg--Soergel consider a category $\Dmix_{(B)}(X)$ of $(B)$-constructible \emph{mixed sheaves} on $X,$ which is a \emph{graded version} of the  $(B)$-constructible derived category $\D_{(B)}(X)$ of sheaves on the flag manifold $X(\C)$ and, equivalently, the bounded derived category of category $\pazocal{O}$ of the Lie algebra $\mathfrak{g}=\operatorname{Lie}(G(\C))$. In particular, there is an autoequivalence $(1)$ of $\Dmix_{(B)}(X)$ called \emph{Tate twist}, shifting the grading, and a functor $$v:\Dmix_{(B)}(X)\to \D_{(B)}(X),$$ called \emph{Betti realisation}, forgetting the grading. 

Most remarkably, Beilinson--Ginzburg--Soergel construct a triangulated equivalence, called \emph{Koszul duality},
$$\gKos: \Dmix_{(B)}(X)\to \Dmix_{(B^\vee)}(X^\vee)$$
mapping projective perverse sheaves to intersection cohomology complexes and intertwining the Tate twist $(1)$ with the shift twist $(1)[2].$ Although very desirable, there is no geometric construction of Koszul duality yet. All six functors commute with Tate twists and hence some new geometric constructions have to enter the picture. \todo{A2}

The main idea of this article is to construct an ungraded version of Koszul duality and to thereby fill in the bottom right corner of the above diagram. This is achieved by a $K$-theoretic point of view and the observation that Betti realisation forgets the Tate twist $(1)$ whereas passage to $K$-theory forgets the shift twist $(1)[2].$ 

To make this observation more precise, we use Soergel--Wendt's very satisfying construction of $\Dmix_{(B)}(X)$ as a full subcategory of the category of \emph{motivic sheaves} $\BM(X/\F_p,\Q),$ see \cite{SoeWe}. We define the category $\KM_{(B^\vee)}(X^\vee)$ as a full subcategory of the category of \emph{$K$-motivic sheaves} $\KM(X^\vee/\F_p,\Q)$ analogously. $K$-motivic sheaves can be thought of as a $K$-theoretic cousin of the category of constructible sheaves computing algebraic $K$-theory instead of Betti cohomology.

There is a functor, which we call \emph{Beilinson realisation},
$$\iota: \Dmix_{(B^\vee)}(X^\vee)\to \KM_{(B^\vee)}(X^\vee),$$
expressing Beilinson's realisation that, rationally, motivic cohomology is a graded refinement of algebraic $K$-theory, see \cite{beilinsonNotesMotivicCohomology1987}. There is a natural isomorphism, called \emph{Bott isomorphism}, $\Q\cong \Q(1)[2]$ in $\KM_{(B^\vee)}(X^\vee)$ and hence $\iota$ forgets the shift by $(1)[2].$

The construction of the ungraded Koszul duality functor 
$$\Kos:\D_{(B)}(X)\to \KM_{(B^\vee)}(X^\vee)$$ is just a copy of Soergel's construction of $\gKos.$ The functor send projective perverse sheaves to so-called intersection $K$-theory complexes, a $K$-theoretic version of intersection cohomology complexes. Both sided admit a combinatorial description in terms of the homotopy category of Soergel modules. The commutativity of the diagram is hence immediate.

We proceed as follows. In the second section we recall Cisinski--D\'eglise's \cite{cisinskiTriangulatedCategoriesMixed2019} construction of $K$-motivic sheaves and motivic sheaves. In the third section we consider the categories $\Dmix_{\pazocal{S}}(X)$ and $\KM_{\pazocal{S}}(X)$ for varieties $X$ with an affine stratification $\pazocal{S}.$ We study natural $t$-structures and weight structures on these categories and recall Soergel's Erweiterungssatz. In the fourth and last section we return to the flag variety.  We recall Soergel's construction of the Koszul duality functor and prove the Theorem mentioned above. In the appendix we collect some useful facts about weight structures and $t$-structures.

\begin{remark}
(1) The case of modular coefficients is work in progress joint with Shane Kelly and building on \cite{EK} and \cite{eberhardtIntegralMotivicSheaves2022}.\\\noindent
(2) The construction of an equivariant (both in the sense of Borel and Bredon) version of $K$-motivic sheaves is work in progress. This will allow to consider an ungraded version of Bezrukavnikov--Vilonen's equivariant/monodromic Koszul duality, see \cite{bezrukavnikovKoszulDualityKacMoody2013}, and open new pathways to Soergel's conjecture on Koszul duality for real groups, see \cite{soergelLanglandsPhilosophyKoszul2001a} and \cite{bezrukavnikovKoszulDualityQuasisplit2021}, as well as the quantum $K$-theoretic geometric Satake, see \cite{cautisQuantumKtheoreticGeometric2015} and\cite{cautisQuantumKtheoreticGeometric2018}

Moreover, the author is considering a motivic Springer correspondence involving the affine Hecke algebra which generalizes \cite{riderFormalityNilpotentCone2013},\cite{riderPerverseSheavesNilpotent2016},\cite{riderFormalityLusztigGeneralized2021}, \cite{eberhardtSpringerMotives2021} and \cite{eberhardtMotivicSpringerTheory2022} and provides a derived version of Lusztig's comparison between the graded affine Hecke algebra and the affine Hecke algebra, see \cite{lusztigAffineHeckeAlgebras1989}.\todo{A3}

\end{remark}
\subsection*{Acknowledgements} We warmly thank Shane Kelly for explaining $K$-motives to us. We are grateful for discussions with Wolfgang Soergel and Geordie Williamson. We thank the referee for helpful comments. This publication was written while the author was a guest at the Max Planck Institute for Mathematics in Bonn.

%% file: Inhalt/motives.tex
Motivic sheaves can be viewed as an amalgamation of the topological notion of sheaves on manifolds and the algebro-geometric concept of motivic cohomology. In this section, we give an overview of the most important properties of motivic sheaves and $K$-motivic sheaves, a lesser known variant. We then discuss the Beilinson realisation functor, concrete constructions and Tate motives. \todo{G1 (General Comment 1)}
\subsection{Background on motivic sheaves}
 To every quasi-projective variety $X$ over a perfect field $k,$ one can associate a $\Q$-linear, tensor-triangulated category of \emph{motivic sheaves} $\BM(X).$ There are various equivalent constructions  and names for the category $\BM(X),$ see \cite[Introduction]{cisinskiTriangulatedCategoriesMixed2019}. For $X=\Spec(k),$ the category of motivic sheaves $\DM(k)$ coincides with Voevodsky's triangulated category of mixed motives over $k$, see \cite{voevodskyTriangulatedCategoriesMotives2000}.
 
 The categories $\BM(X)$ are equipped with a six functor formalism admitting, for example, localisation triangles, projection formulae and base change, see \cite[A.5]{cisinskiTriangulatedCategoriesMixed2019}. Using the six functors, one can define the \emph{motive} (with compact support) of a variety $f:X\to\Spec(k)$ as  $$\operatorname{M}(X)=f_*f^!\Q\in \BM(k) \text{ and } \operatorname{M}^c(X)=f_!f^!\Q\in \BM(k)\text{, respectively.}$$
 The motive of the projective line $\mathbb{P}^1_k$ splits as a direct sum
 $$\operatorname{M}(\mathbb{P}^1_k)=\Q\oplus \Q(1)[2]$$
 and induces an autoequivalence $(1)=-\otimes \Q(1)$ called \emph{Tate twist} that commutes with all six functors.
 
 Homomorphisms between motivic sheaves are governed by algebraic cycles. For $X$ smooth, there is a natural isomorphism
$$\Hom{\BM(X,\Q)}{\Q}{\Q(p)[q]}\cong \operatorname{CH}^p(X,2p-q)_\Q$$
 with Bloch's higher Chow groups, see \cite{Blo86}. For $q=2p,$ these are the usual Chow groups $\operatorname{CH}^p(X)_\Q$ of codimension-$p$ algebraic cycles on $X$ up to rational equivalence. For cellular varieties, such as flag varieties, considered in this article these Chow groups will just coincide with the usual Borel--Moore homology of their complex points.
 
For each prime $\ell$ invertible in $k$, there is an \emph{$\ell$-adic realisation functor}
$$\operatorname{Real}_\ell:\BM(X,\Q)\to \operatorname{D}_{\text{\'et}}(X, \Q_\ell)$$
to the category of $\ell$-adic sheaves and for $k=\C$ a \emph{Betti realisation functor}
 $$v:\BM(X,\Q)\to \operatorname{D}(X(\C), \Q)$$
 to the category of sheaves on $X(\C)$ equipped with the metric topology, see \cite{ayoubRealisationEtaleOperations2014} and \cite{ayoubNoteOperationsGrothendieck2009}. Both types of realisation functors are compatible with the six functors and induce the cycle class maps
 $$\operatorname{CH}^p(X)_\Q\to H^{2p}_{\text{\'et}}(X,\Q_\ell) \text{ and } \operatorname{CH}^p(X)_\Q\to H^{2p}_{\text{Betti}}(X(\C),\Q)\text{, respectively}.$$
For our purposes, all varieties are defined over $\Z$ and all sheaves we consider are of geometric origin. We will hence take the freedom to treat $\ell$-adic sheaves on $X/\overline{\mathbb{F}}_p$ and sheaves on $X(\mathbb{C})$ interchangeably using \cite{BBD}.
\subsection{Background on $K$-motivic sheaves} Very similarly to the system of categories of motivic sheaves $\DM(X)$ there is a system of \emph{$K$-motivic sheaves} $\KM(X)$ associated to quasi-projective varieties $X$ over a perfect field $k,$ see \cite[Section 13.3]{cisinskiTriangulatedCategoriesMixed2019}. $K$-motivic sheaves are also equipped with a full six-functor formalism and Tate twist with the same properties. 

A main difference is that the $K$-motive of $f:\mathbb{P}^1_k\to\Spec(k)$ splits into two isomorphic copies
$$\operatorname{M}(\mathbb{P}^1_k)=f_*f^!\Q=\Q\oplus \Q\in \KM(k).$$
This implies that the shift twist functor $(1)[2]$ acts as the identity on $\KM(X),$ a phenomenon known as \emph{Bott periodicity.} 

As the name suggests, homomorphisms between $K$-motivic sheaves are governed by $K$-theory. For $X$ smooth, there is natural isomorphism 
$$\Hom{\KM(S,\Q)}{\Q}{\Q(p)[q]}=K_{2p-q}(X)_\Q$$
with the rational higher algebraic $K$-theory of $X.$
\subsection{The Beilinson realisation functor} There is a very close relationship between motivic cohomology (which is isomorphic to higher Chow groups) and algebraic $K$-theory, first observed by Beilinson \cite{beilinsonNotesMotivicCohomology1987}. 
The rational $K$-theory of a smooth variety $X$ naturally decomposes in eigenspaces of the Adams operations
\todo{S3}
\begin{equation}\label{eq:decompktheory}
K_n(X)_\Q=\bigoplus_{i} K^{(i)}_n(X),
\end{equation}
turning $K_\bullet(X)$ into a bigraded ring, see \cite[IV.5]{Wei}. By \cite{Blo86}, there are natural isomorphisms
$K^{(i)}_{n}(X)\cong \operatorname{CH}^i(X,n)_\Q$ such that for $n=0$ the decomposition in \eqref{eq:decompktheory} yields the Chern character isomorphism
$$K_0(X)_\Q\cong \bigoplus_i \operatorname{CH}^i(X)_\Q.$$
Hence, rational motivic cohomology can be regarded as a graded refinement of algebraic $K$-theory. These results admit a relative version. Namely, there is functor $$\iota: \BM(X)\to \KM(X),$$
which we call \emph{Beilinson realisation}. It is compatible with the six functors and Tate twists and it induces for all $M,N\in\DM(X)$ an isomorphism
\begin{align}\label{eq:iotagrading}
\Hom{\KM(X)}{\iota(M)}{\iota(N)}&\cong\bigoplus_i\Hom{\BM(X)}{M}{N(i)[2i]}
\end{align}
that specialises to the decomposition \eqref{eq:decompktheory} in Adams eigenspaces. This way, motivic sheaves can be regarded as a graded refinement of $K$-motivic sheaves where $(1)[2]$ is the shift of grading functor and $\iota$ forgets the grading.
\label{sec:beilisonmotivesandkmotives}
\subsection{Construction of ($K$-)motivic sheaves}
We sketch Cisinski--D\'eglise's construction of the categories of $K$-motivic sheaves $\KM(X)$ and motivic sheaves $\BM(X),$  \cite[Chapter 13-15]{cisinskiTriangulatedCategoriesMixed2019}.

First, one considers the ring spectrum $KGL_{\Q,X}$ representing rational homotopy invariant $K$-theory in the stable motivic homotopy category $\operatorname{SH(X)}.$ This allows to consider the category of \emph{$K$-motivic sheaves} over $X,$
$$\KM(X):= \operatorname{Ho}(KGL_{\Q,X}\modules),$$
as the homotopy category of modules over $KGL_{\Q,X}.$ The system of categories $\KM(X)$ forms a so-called \emph{motivic triangulated category} (see \cite[Definition 2.4.45]{cisinskiTriangulatedCategoriesMixed2019}) which entails a six functor formalism with all desired properties.



As shown by Riou \cite{Riou}, the spectrum $KGL_X$ admits an Adam's decomposition similar to \eqref{eq:decompktheory}$$KGL_{\Q,X}=\bigoplus_{i}KGL^{(i)}_X.$$
Denoting $H_{\text{\hspace{-10pt}\cyr B},X}:=KGL^{(0)}_X,$  Cisinski--D\'eglise define the category of \emph{Beilinson motives} over $X,$
$$\BM_\text{\hspace{-10pt}\cyr B}(X):=\operatorname{Ho}(H_{\text{\hspace{-10pt}\cyr B},X}\modules),$$
as the homotopy category of modules over $H_{\text{\hspace{-10pt}\cyr B},X}.$ The system of categories $\BM_\text{\hspace{-10pt}\cyr B}(X)$ forms a motivic triangulated category and is shown to be equivalent to other definitions of motivic sheaves, see \cite[Chapter 16]{cisinskiTriangulatedCategoriesMixed2019}. We hence write $\BM(X)=\BM_\text{\hspace{-10pt}\cyr B}(X)$ and refer to objects as motivic sheaves.

The inclusion map $H_{\text{\hspace{-10pt}\cyr B},X}\to KGL_{\Q,S}$ induces a forgetful functor $F:\KM(X)\to \BM(X)$ whose left adjoint \todo{S4}
$$\iota: \BM(X,\Q)\rightarrow \KM(S,\Q),\, M\mapsto KGL_{\Q,X}\otimes_{H_{B,X}}M$$
we call  \emph{Beilinson realisation functor.} 
By \cite[Lemma 14.1.4]{cisinskiTriangulatedCategoriesMixed2019} each $KGL_{S,\Q}^{(i)}$ is naturally isomorphic to $H_{\text{\hspace{-10pt}\cyr B},S}$ as an $H_{\text{\hspace{-10pt}\cyr B},S}$-module and there is an isomorphism of ring spectra
$$KGL_{\Q,X}\cong H_{\text{\hspace{-10pt}\cyr B},X}[t,t^{-1}]:=\bigoplus_{i\in\Z} H_{\text{\hspace{-10pt}\cyr B},X}(i)[2i].$$
For $M\in \BM(X)$ this yields the following simple formula \todo{S4}
$$F\iota(M)\cong\bigoplus_{i\in\Z} M(i)[2i]$$
which implies \eqref{eq:iotagrading}.
Using \cite[Theorem 4.4.25]{cisinskiTriangulatedCategoriesMixed2019} one shows that $\iota$ is compatible with the six functors. For further properties on $K$-motivic sheaves we refer to \cite[Proposition 4.1.1]{BL}.
\subsection{Tate Motives over Affine Spaces of Finite Fields.}\label{sec:Tateoverfield}Let $k=\F_q$ be a finite field. We denote by 
\begin{align*}
\BMT(\A^n_k, \Q)&=\left\langle \Q(n)\,|\, n\in \Z\right\rangle_\Delta\subset \BM(\A^n_k,\Q) \text{ and }\\
\KMT(\A^n_k, \Q)&=\left\langle \Q \right\rangle_\Delta\subset \KM(\A^n_k,\Q)
\end{align*}
the categories of \emph{mixed Tate ($K$-)motives} (observe that $\Q(n)\cong \Q[-2n]$ in  $\KM$ and is hence not needed as a generator). 
Since the rational higher $K$-theory and the rational higher Chow groups of a finite field vanish, these categories of become semi-simple and one can easily show:
\begin{theorem} There are equivalences of tensor-triangulated categories
\begin{align*}
\BMT(\A^n_k, \Q)&\cong\Der^b(\Q\modules^\Z)\text{ and}\\
\KMT(\A^n_k, \Q)&\cong\Der^b(\Q\modules)
\end{align*}
with the bounded derived categories of (graded) finite dimensional vector spaces over $\Q.$ Here we let $\Q(p)\in \BMT(\A^n_k, \Q)$ correspond to $\Q$ sitting in grading degree $-p$ and cohomological degree $0$ by convention.
\end{theorem}
This equips the categories $\BMT(\A^n_k, \Q)$ and $\KMT(\A^n_k, \Q)$ with canonical $t$-structures we denote by $t$, see Remark \ref{rem:tstructuresonderivedcategories}(1).

For quasi-projective  varieties $X/k$  the categories $\BM(X, \Q)$ and $\KM(X, \Q)$ are equipped with a weight structure $w,$ see \cite{Heb} and \cite{BL}.  This weight structure descends to Tate motives and assigns the weight $2p-q$ to $\Q(p)[q]$ such that
$$\Q(p)[q]\in \BM(X)^{w=2p-q} \text{ and } \Q(p)[q]\in \KM(X)^{w=2p-q}.$$

We observe that the $t$-structure and weight structure on $\KMT(\A^n_k)$ coincide. As explained in Remark \ref{rem:tstructuresonderivedcategories}(1), $t$-structures and weight structures usually behave very differently. Our case just happens to be quite degenerate, since $\KMT(\A^n_k)$ is semi-simple.

We observe that the induced functor
$$\iota: \BMT(\A^n_k)\rightarrow \KMT(\A^n_k)$$
is \emph{not} compatible with the $t$-structures, since $\iota(\Q(1)[2])=\Q.$ Rather, $\iota$ is compatible with the weight structures:
\begin{proposition}\label{prop:weighttonpoint} Let $M$ be in $\KMT(\A^n_k)$, then
\begin{align*}
M\in \BMT(\A^n_k)^{w\leq 0}& \iff \iota(M)\in  \KMT(\A^n_k)^{w\leq 0}=\KMT(\A^n_k)^{t\leq 0} \text{ and}\\
M\in \BMT(\A^n_k)^{w\geq 0}& \iff \iota(M)\in \KMT(\A^n_k)^{w\geq 0}=\KMT(\A^n_k)^{t\geq 0}.
\end{align*}
\end{proposition}
As explained in \cite[Section 3.4]{SoeWe}, the interplay of the $t$-structure and weight structure on $\BMT(\A^n_k, \Q)$ can be seen as toy case of Koszul duality. So Proposition \ref{prop:weighttonpoint}  gives a subtle first hint that $\iota$ should be related to Koszul duality!

%% file: Inhalt/affinelystratified.tex
\subsection{Constructible motivic sheaves} Let $k=\F_q$ be a finite field. Let $X/k$ be a variety with a cell decomposition (also called affine stratification), that is, $$X=\biguplus_{s\in \Ss} X_s$$ 
where $\Ss$ is some finite set and each $i_s:X_s\to X$ is a locally closed subvariety isomorphic to $\A^n_k$ for some $n\geq 0.$ In this situation, Soergel--Wendt \cite{SoeWe} make the following definition: 
\begin{definition}The category of \emph{mixed stratified Tate motives} on $X$ is
$$\MTDer{\Ss}{X,\Q}=\left\{M\in \BM(X,\Q) \,|\, i_s^*M\in \BMT(X_s,\Q) \text{ for all } s\in\Ss\right\}$$
the full subcategory of motivic sheaves which restrict to mixed Tate motives on the strata.
\end{definition}
For this category to be well-behaved, so for example closed under Verdier duality,  Soergel--Wendt impose the following technical condition on the stratification.
\begin{definition} The stratification $\Ss$ on $X$ is called \emph{Whitney--Tate} if $i_t^{*}i_{s,*}\Q\in \BMT(X_t,\Q)$ for all $s,t\in\Ss.$
\end{definition}
We will abbreviate $\Dmix_\Ss(X)=\MTDer{\Ss}{X,\Q}$ 
and speak of \emph{$\Ss$-constructible motivic sheaves,} 
and assume that $\Ss$ is Whitney--Tate from now on.

We can now copy their definition in the context of $K$-motives.
\begin{definition} The category of \emph{$\Ss$-constructible $K$-motivic sheaves} is
$$\KM_\Ss(X)=\left\{M\in \KM(X,\Q) \,|\, i_s^*M\in \KMT(X_s,\Q) \text{ for all } s\in\Ss\right.\}$$
the full subcategory of the category of $K$-motivic sheaves $\KM(X,\Q)$ of objects which restrict to Tate motives on the strata.
\end{definition}

Since the functor $\iota:\BM(X)\rightarrow \KM(X)$ commutes with the six operations, we see that it descends to a functor
$$\iota: \Dmix_\Ss(X)\to \KM_\Ss(X)$$
and observe that the Whitney--Tate condition with respect to $\BM(X)$ implies the one for $\KM(X).$

In order to be closed under the six functors, we need to restrict us to morphisms of varieties which are compatible with their affine stratification in the following sense.
\begin{definition}
	Let $(X,\Ss)$ and $(Y,\Ss^\prime)$ be varieties with affine stratifications. We call
	$f:X\rightarrow Y$ an \emph{affinely stratified map} if
	\begin{enumerate}
		\item for all $s \in\Ss^\prime$ the inverse image $f^{-1}(Y_s)$ is a union of strata;
		\item for each $X_s$ mapping into $Y_{s^\prime}$, the induced map $f:X_s\rightarrow Y_{s^\prime}$ is a surjective affine map. \todo{S5}
	\end{enumerate}
\end{definition}
\subsection{Weight Structures}
The categories $\DM_{gm}(X)$ and $\KM_{gm}(X)$ of objects of geometric origin naturally come with a weight structure, called \emph{Chow weight structure}, whose hearts are generated by objects of the form $f_*\Q$ for smooth projective maps $f:Y\to X,$ see \cite{Heb}.

We will still define the weight structures on the subcategories $\Dmix_\Ss(X)$ and $\KM_\Ss(X)$ by hand using the gluing formalism described in Appendix \ref{sec:tandweightstructures}. We note that our definition coincides with the restriction of the Chow weight structures. \todo{S6}
\begin{theorem} 
Setting
\begin{align*}
\Dmix_\Ss(X)^{w\leq 0}&=\left\{M\in \Dmix_\Ss(X) \,|\, i_s^!M\in \BMT(X_s)^{w\leq 0}\text{ for all } s\in \Ss\right\},\\
\Dmix_\Ss(X)^{w\geq 0}&=\left\{M\in \Dmix_\Ss(X) \,|\, i_s^*M\in\BMT(X_s)^{w\geq 0}\text{ for all } s\in \Ss\right\},\\
\KM_\Ss(X)^{w\leq 0}&=\left\{M\in \KM_\Ss(X) \,|\, i_s^!M\in \KMT(X_s)^{w\leq 0}\text{ for all } s\in \Ss\right\}\text{ and }\\
\KM_\Ss(X)^{w\geq 0}&=\left\{M\in \KM_\Ss(X) \,|\, i_s^*M\in\KMT(X_s)^{w\geq 0}\text{ for all } s\in \Ss\right\}
\end{align*}
defines bounded weight structures on $\Dmix_\Ss(X)$ and $\KM_\Ss(X)$. \todo{A5}
\end{theorem}
\begin{proof} Use Theorem \ref{thm:glueing} inductively.
\end{proof}
The weight structures on $\Dmix_\Ss(X)$ and $\KM_\Ss(X)$ are closely related:
\begin{proposition}\label{prop:weightt} Let $M\in \Dmix_\Ss(X).$ Then 
\begin{align*}
M\in \Dmix_\Ss(X)^{w\leq 0}& \text{ if and only if } \iota(M)\in \KM_\Ss(X)^{w\leq 0} \text{ and}\\
M\in \Dmix_\Ss(X)^{w\geq 0}& \text{ if and only if } \iota(M)\in \KM_\Ss(X)^{w\geq 0}.
\end{align*}
\end{proposition}
\begin{proof}
Since $\iota$ commutes with $i_s^*$ and $i_s^!$ the statement follows from Proposition \ref{prop:weighttonpoint}. 
\end{proof}
Moreover, the six functors have certain exactness properties with respect to weight structures:
\begin{theorem}\label{thm:exactness} Let $(X,\Ss)$ and $(Y,\Ss')$ be varieties with a Whitney--Tate affine stratification and $f:X\rightarrow Y$ be an affinely stratified map. \todo{S8} Then we get
\begin{enumerate}
\item $f^*, f_!$ and $\otimes$ are right $w$-exact.
\item $f^!, f_*$ are left $w$-exact.
\end{enumerate}
\end{theorem}
\begin{proof}
See \cite[Proposition 3.2]{EK}.
\end{proof}
\subsection{Pointwise purity and the weight complex functor}
In general, objects in the heart of a weight structure $w$ on a category $\Cc$ admit no positive extension, that is,  $\Hom{\Cc}{M}{N[n]}=0$ for all $M,N\in \Cc^{w=0}$ and $n>0.$ We will show that under a certain \emph{pointwise purity} assumption there are also no negative extensions for objects in $\Dmix_\Ss(X)^{w=0}$ and $\KM_\Ss(X)^{w=0}.$ \todo{A4}
\begin{definition} Let $M$ be in $\Dmix_\Ss(X)$ (resp. $\KM_\Ss(X)$). We say that $M$ is \emph{pointwise pure} if
$i_s^*M$ and $i_s^!M$ are in $\BMT(X_s)^{w=0}$ (resp.  $\KMT(X_s)^{w=0}$) for all $s\in\Ss.$
\end{definition}
\begin{proposition} $M\in \Dmix_\Ss(X)$ is pointwise pure if and only if $\iota(M)\in \KM_\Ss(X)$ is pointwise pure.
\end{proposition}
One can construct pointwise pure objects by different methods, for example:
\begin{enumerate}
\item Using affinely stratified resolutions of singularities of closures of strata in $X,$ see \cite[Theorem 4.5]{EK}.
\item Using contracting $\mathbb{G}_m$ actions, see \cite[Proposition 7.3.]{SoeWe}.
\item By an inductive process in the case of flag varieties, see \cite[Lemma 6.6]{SoeWe}.
\end{enumerate}
All of those methods can be used to show that all objects in $\Dmix_\Ss(X)^{w=0}$ and $\KM_\Ss(X)^{w=0}$ are pointwise pure in the case of flag varieties.

Pointwise pure objects in $\Dmix_\Ss(X)$ are under some assumptions in fact sums of (appropriately \todo{S7} shifted and twisted) intersection complexes, that is, simple perverse motives. See \cite[Corollary 11.11]{SoeWe}. 

Pointwise pure objects are very special since they have no non-trivial extensions amongst each other.
\begin{proposition}\label{prop:pointwisepuritynoexts} Let $M,N$ be in $\Dmix_\Ss(X)$ (resp. $\KM_\Ss(X)$) be pointwise pure. Then for all $n\neq 0$ we have
$$\Hom{\Dmix_\Ss(X)}{M}{N[n]}=0 \text{ (resp. }\Hom{\KM_\Ss(X)}{M}{N[n]}=0\text{).}$$
\end{proposition}
\begin{proof} The statement for $\Dmix_\Ss(X)$ follows from the one of $\KM_\Ss(X)$ using $\iota.$ For $\KM_\Ss(X),$ we observe that the pointwise purity implies that $M,N\in \KM_\Ss(X)^{w=0}\cap \KM_\Ss(X)^{t=0},$ where by $\KM_\Ss(X)^{t=0}$ we denote the heart of the bottom $p=0$ perverse $t$-structure on $\KM_\Ss(X).$ Hence the statement for negative $n$ follows from the axioms of the $t$-structure and the statement for positive $n$ from the axioms of the weight structure.
\end{proof}
Pointwise purity allows us to consider the category $\Dmix_\Ss(X)$ and $\KM_\Ss(X)$ as homotopy categories of their weight zero objects. 
\begin{theorem}\label{thm:weightcomplexDmix} Assume that all objects in  $\KM_\Ss(X)^{w=0}$ are pointwise pointwise pure. Then the \emph{weight complex functor}(see Theorem \ref{thm:weightcomplex}) induces an equivalences of categories,
\begin{align*}
\Dmix_\Ss(X)&\cong\Hot^b(\Dmix_\Ss(X)^{w=0})\\
\KM_\Ss(X)&\cong\Hot^b(\KM_\Ss(X)^{w=0})
\end{align*}
compatible with the functor $\iota.$
\end{theorem}
\begin{proof}
We prove the statement for $\Dmix_\Ss(X),$ the case of $\KM_\Ss(X)$ is done in the same way. The pointwise purity assumption and Proposition \ref{prop:pointwisepuritynoexts} shows that there are no non-trivial extensions in $\Dmix_\Ss(X)$ between objects in $\Dmix_\Ss(X)^{w=0}.$ Trivially, the same holds true in $\Hot^b(\Dmix_\Ss(X)^{w=0}).$ Since the weight complex functor restricts to the inclusion $\Dmix_\Ss(X)^{w=0}\rightarrow \Hot(\Dmix_\Ss(X))$ an inductive argument (``d\'evissage'') shows that the functor is indeed fully faithful, where we use that $\Dmix_\Ss(X)$ is generated by $\Dmix_\Ss(X)^{w=0}$ as a triangulated category. \todo{A5}

The compatibility with $\iota$ follows since $\iota$ is weight exact.
\end{proof}
We note that there is a different way of proving the last theorem using a formalism called ``tilting'', see \cite{SoeWe} and \cite{SVW}. We prefer the weight complex functor, since it also exists without the pointwise purity assumption. The weight complex functor even exists for all motivic sheaves and $K$-motivic sheaves of geometric origin, where the heart of the weight structure is the category Chow motives, see \cite{Bon}.

\subsection{Erweiterungssatz} \label{sec:erweiterungssatz}
The \emph{Erweiterungssatz} as first stated in \cite{Soe90} and reproven in a more general setting in \cite{Gin} allows a \emph{combinatorial} description of pointwise pure weight zero sheaves on $X$ in terms of certain modules over the cohomology ring of $X$. In the case of $X$ being the flag variety, these modules are called \emph{Soergel modules}. In \cite{SoeWe} a motivic version is considered, which easily extends to $K$-motives.
\begin{definition} We denote by 
	\begin{align*}
\Hyp: \Dmix_\Ss(X)\rightarrow\Hyp(X)\modules^{\Z},\, &M\mapsto\bigoplus_{n\in\Z}\Hom{\Dmix_\Ss(X)}{\Q}{M(n)[2n]}\\
\Kyp: \KM_\Ss(X)\rightarrow\Kyp(X)\modules,\, &M\mapsto\Hom{\KM_\Ss(X)}{\Q}{M}
\end{align*}
the \emph{hypercohomology} functors. 
Here $\Hyp(X)=\bigoplus_{n\in\Z}\Hom{\Dmix_\Ss(X)}{\Q}{\Q(n)[2n])},$ and $\Kyp(X)=\Hom{\KM_\Ss(X)}{\Q}{\Q},$ and the former is interpreted as a graded ring.
\end{definition}
The rings $\Hyp(X)$ and $\Kyp(X)$ are nothing else than the motivic cohomology and $K$-theory of $X$ and we collect some of the important properties.
\begin{theorem}
(1) The map $\iota: \Hyp(X)\rightarrow \Kyp(X)$ induced by $\iota: \Dmix_\Ss(X)\to \KM_\Ss(X)$ is an isomorphism.\\\noindent
(2) The following diagram commutes up to natural transformation
\begin{center}\begin{tikzcd}
\Dmix_\Ss(X) \arrow{r}{\Hyp} \arrow{d}{\iota}
& \Hyp(X)\modules^{\Z} \arrow{d} \\
\KM_\Ss(X)  \arrow{r}{\Kyp}
& \Kyp(X)\modules
\end{tikzcd}\end{center}
where the right vertical arrow is forgetting the grading.\\\noindent
(3) The ring $\Hyp(X)$ is the Chow ring of $X$. Assume that $X$ and the stratification is already defined over $\Z.$ Then $\Hyp(X)$ coincides with the Borel--Moore singular homology of $X^{an}(\C).$\\\noindent
(4) Assume that $X$ is smooth, then $\Kyp(X)=K_0(X)$ is $0$-th K-group of $X,$ that is, the Grothendieck group of the category of vector bundles on $X.$
\end{theorem}
\begin{proof}All statements are direct consequences of the discussion in Section \ref{sec:beilisonmotivesandkmotives}.
\end{proof}
We remark that motivic cohomology is bigraded (higher Chow groups) and $K$-theory graded (higher $K$-groups). In our particular setup (affine stratification, finite field base, rational coefficients) all the higher groups vanish. We hence see one grading less.

Under a certain technical assumption the functors $\Hyp$ and $\Kyp$ are fully faithful on pointwise pure objects.
\begin{theorem}[Erweiterungssatz]\label{thm:erweiterungssatz} Assume that all objects in  $\Dmix_\Ss(X)^{w=0}$ are pointwise pure and for each stratum $i:X_s\to X$  and $M\in  \Dmix_\Ss(X)^{w=0}$ the map $\Hyp(M)\to \Hyp(i_*i^*M)$ is surjective and the map $\Hyp(i_!i^!M)\to\Hyp(M)$ is injective.
Then the functors
\begin{align*}
\Hyp&: \Dmix_\Ss(X)^{w=0}\rightarrow\Hyp(X)\modules^{\Z}\\
\Kyp&: \KM_\Ss(X)^{w=0}\rightarrow\Kyp(X)\modules
\end{align*}
are fully faithful.
\end{theorem}
\begin{proof}
The statement for $\Dmix_\Ss(X)$ is proven in \cite[Section 8]{SoeWe}. The proof uses the six functor formalism and weight arguments and applies word for word for $\KM_\Ss(X)^{w=0}.$
\end{proof}
The assumptions are fulfilled if there are contracting $\mathbb{G}_m$ actions for the closure of strata, see \cite[Proposition 8.8]{SoeWe}. The theorem in particular applies to flag varieties.
\begin{definition} We denote the essential images of $\Dmix_\Ss(X)^{w=0}$ and $\KM_\Ss(X)^{w=0}$  under $\Hyp$ and $\Kyp$ by $\Hyp(X)\Smodules^\Z$ and  $\Kyp(X)\Smodules^\Z$ and call them \emph{Soergel modules.}
\end{definition}
\begin{corollary} Under the assumptions of Theorem \ref{thm:erweiterungssatz}, there are equivalences of categories
\begin{align*}
\Dmix_\Ss(X)&\cong \Hot^b(\Hyp(X)\Smodules^\Z)\text{ and}\\
\KM_\Ss(X)&\cong  \Hot^b(\Kyp(X)\Smodules).
\end{align*}
\end{corollary}

%% file: Inhalt/flagvariety.tex
We discuss the particular case of flag varieties and Koszul duality.
\subsection{Flag Varieties} Let $G\supset B \supset T$ be a split reductive group over $\F_p$ with a Borel subgroup and maximal torus. Denote the Langlands dual by $G^\vee\supset B^\vee \supset T^\vee.$ Denote by $X(T)=\Hom{}{T}{\Gm}$ the character lattice, by $W=N_G(T)/T\supset S$ the Weyl group with the set of simple reflections corresponding to $B$, and for $w\in W$ by $l(w)\in \Z_{\geq 0}$ the length of an element. The flag variety $X=G/B$ has an affine stratification $(B)$ by its $B$-orbits, called the \emph{Bruhat stratification}, 
$$X=\biguplus_{w\in W} X_w$$
where $X_w=BwB/B\cong\A_k^{l(w)}.$ By \cite[Proposition 4.10]{SoeWe} this stratification fulfills the Whitney--Tate condition. More generally, the partial flag varieties $G/P$ for parabolics $B\subset P \subset G$ with their stratification by $B$-orbits are Whitney--Tate.
Similarly, $X^\vee=G^\vee/B^\vee$ admits a Bruhat stratification by $B^\vee$-orbits.
It hence makes sense to consider the categories $\Dmix_{(B)}(X),$  $\D_{(B)}(X)$ as well as $\Dmix_{(B^\vee)}(X^\vee),$  $\KM_{(B^\vee)}(X^\vee)$ from the Introduction.

\subsection{Translation Functors and Pointwise Purity} We recall the inductive construction of pointwise pure objects in $\Dmix_{(B^\vee)}(X^\vee),$ see \cite[Section 6]{SoeWe}. 

First of all, the object $i_{e,!}\Q$ is pointwise pure, where $e\in W$ denotes the identity. For a simple reflection $s\in S$ we denote by $P^\vee_s=B^\vee\cup B^\vee s B^\vee$ the minmal parabolic and the smooth proper morphism (in fact the map is a projective bundle)
$$\pi_s: X^\vee\to G^\vee/P^\vee_s.$$
The functor $\theta_s=\pi_s^*\pi_{s,*}$ is called \emph{translation functor}. It clearly preserves pointwise pure objects. For an arbitary $w\in W$ with $l=l(w)$ we choose a shortest expression $w=s_1\cdots s_l.$ Then the object
$\theta_{s_1}\cdots\theta_{s_n}i_{e,!}\Q$ is called a \emph{Bott--Samelson motive}. It is pointwise pure, has support $\overline{X^\vee_w}$ and a unique indecomposable direct summand, which we will denote by $\gEc_w,$ with support $\overline{X^\vee_w}.$ In fact all pointwise pure objects are sums of shifts twits of the motives $\gEc_w$ and the objects $\mathcal{IC}_w=\gEc_w[l(w)]$ are \emph{simple perverse motives (intersection  cohomology complexes)} by the decomposition theorem! Subsumed, we get
\begin{align*}
\Dmix_{(B^\vee)}(X^\vee)^{w=0}&=\genbuild{\theta_{s_1}\cdots\theta_{s_l}i_{e,!}\Q(n)[2n]}{s_i\in S, l\in \Z_{\geq 0}, n\in \Z}_{\cong, \oplus, \inplus}\\
&=\genbuild{\gEc_w(n)[2n]}{w\in W, n\in \Z}_{\cong, \oplus}.
\end{align*}
where by $\cong, \oplus, \inplus$ we denote closure under isomorphism, finite direct sum and direct summands.
We see that all objects in $\Dmix_{(B^\vee)}(X^\vee)^{w=0}$ are pointwise pure.

We observe that exactly the same construction works for $\KM_{(B^\vee)}(X^\vee).$ We denote $\Ec_w=\iota( \gEc_w)$ and obtain
\begin{align*}
\KM_{(B^\vee)}(X^\vee)^{w=0}&=\genbuild{\theta_{s_1}\cdots\theta_{s_l}i_{e,!}\Q}{s_i\in S, l\in \Z_{\geq 0}}_{\cong, \oplus, \inplus}\\
&=\genbuild{\Ec_w}{w\in W}_{\cong, \oplus}.
\end{align*}
Here, the objects $\Ec_w$ may be taken as a definition for \emph{intersection $K$-theory complexes} $\mathcal{IK}_w=\Ec_w.$

Furthermore, the weight complex functor induces equivalences of categories, see Theorem \ref{thm:weightcomplexDmix},
\begin{align*}
\Dmix_{(B^\vee)}(X^\vee)\cong \Hot^b(\Dmix_{(B^\vee)}(X^\vee)^{w=0})&\text{ and}\\
\KM_{(B^\vee)}(X^\vee)\cong \Hot^b(\KM_{(B^\vee)}(X^\vee)^{w=0})
\end{align*}
compatible with the functor $\iota$ in the obvious way.
\subsection{Soergel Modules I} The categories $\Dmix_{(B^\vee)}(X^\vee)^{w=0}$ and $\KM_{(B^\vee)}(X^\vee)^{w=0}$ can be described combinatorially in terms of Soergel modules, using the functors $\Hyp$ and $\Kyp$ and the Erweiterungssatz, see Section \ref{sec:erweiterungssatz}. We recall the explicit description of $\Hyp(X^\vee)$ and $\Kyp(X^\vee).$ Recall that $X(T)=\Hom{}{T^\vee}{\Gm}$ denotes the character lattice. Then there are natural isomorphisms
$$C=\Sym(X(T^\vee)\otimes\Q)/\Sym(X(T^\vee)\otimes\Q)^W_+\cong\Hyp(X^\vee)\cong\Kyp(X^\vee)$$
where $\Sym(X(T)\otimes\Q)$ denotes the symmetric algebra, and $\Sym(X(T^\vee)\otimes\Q)^W_+$ the ideal of invariants of positive degree under the action of $W$. In \cite{Soe90} it is shown that the Bott--Samelson motives
$$\theta_{s_1}\cdots\theta_{s_l}i_{e,!}\Q$$
are mapped to the Bott--Samelson modules
$$C\otimes_{C^{s_1}}\cdots\otimes_{C^{s_l}}\Q$$
under $\Hyp$ and hence also under $\Kyp.$ We denote $D_w=\Hyp(\gEc_w)=\Kyp(\Ec_w).$ Then $D_w$ can be characterised as the unique indecomposable direct summand of the Bott-Samelson module not appearing in the Bott-Samelson modules associated to elements of $W$ of shorter length. 

Furthermore, the assumption for the Erweiterungssatz (Theorem \ref{thm:erweiterungssatz}) are fulfilled by flag varieties, see \cite[Proposition 8.8]{SoeWe}. Hence, there is an equivalence of categories
\begin{align*}
\Dmix_{(B^\vee)}(X^\vee)^{w=0}&=\genbuild{\gEc_w(n)[2n]}{w\in W, n\in \Z}_{\cong, \oplus}\\
&\cong\genbuild{D_w\langle n\rangle}{w\in W, n\in \Z}_{\cong, \oplus}\\
&=C\Smodules^\Z
\end{align*}
of weight zero objects in $\Dmix_{(B^\vee)}(X^\vee)$ and the category of graded \emph{Soergel modules}. Here we denote by $\langle n\rangle$ the shift of grading in $C\modules^\Z.$  In the same way there is an equivalence
\begin{align*}
\KM_{(B^\vee)}(X^\vee)^{w=0}&=\genbuild{\Ec_w}{w\in W,}_{\cong, \oplus}\\
&\cong\genbuild{D_w}{w\in W}_{\cong, \oplus}\\
&=C\Smodules
\end{align*}
between pointwise pure $K$-motives and ungraded \emph{Soergel modules}. Both descriptions are compatible with the functor $\iota$ in the obvious way.
\subsection{Projective Perverse Sheaves} We describe the ``Koszul dual'' of the last sections. This is the ``classical story''. First, there is a functor, called \emph{Betti realization functor}, 
$$v: \Dmix_{(B)}(X)\rightarrow \D_{(B)}(X)$$
where $\D_{(B)}(X)=\D_{(B)}(X(\C))$ is the $(B)$-constructible derived category of sheaves on the complex manifold $X(\C)$. The functor $v$ is a degrading functor with respect to the Tate twist $(n)$. We have $v(-(n))\cong v(-)$ and for $M, N \in \Dmix_{(B)}(X)$ the functor $v$ induces an isomorphism
$$\bigoplus_{n\in \Z}\Hom{\Dmix_{(B)}(X)}{M}{N(n)}=\Hom{\D_{(B)}(X)}{v(M)}{v(N)}.$$
The functor $v$ is furthermore clearly exact for the perverse $t$-structures on $\Dmix_{(B)}(X)$ and $\D_{(B)}(X).$ 
We denote the categories of projective perverse sheaves by $$\Proje \Dmix_{(B)}(X)^{t=0}\text{ and }\Proje \D_{(B)}(X)^{t=0}.$$
One can show, using Theorem \ref{thm:weightcomplex}, that there are equivalences of categories
\begin{align*}
\Hot^b(\Proje \Dmix_{(B)}(X)^{t=0})&\cong\Der^b(\Dmix_{(B)}(X)^{t=0})\cong \Dmix_{(B)}(X)\text{ and}\\
\Hot^b(\Proje \D_{(B)}(X)^{t=0})&\cong\Der^b(\D_{(B)}(X)^{t=0})\cong \D_{(B)}(X)
\end{align*}
which are all compatible with $v.$
Denote by $w_0\in W$ the longest element. Let $\gPc_w\in \Proje \Dmix_{(B)}(X)^{t=0}$ be the projective cover of $j_{ww_0,!}\Q[l(w)].$ Then $\Pc_w=v(\gPc_w)\in \Proje \D_{(B)}(X)^{t=0}$ the projective cover of $j_{ww_0,!}\Q[l(w)]\in D(X)^{t=0}.$
\subsection{Soergel Modules II}
Soergel shows in \cite{Soe90} that categories of projective perverse objects  $\Proje \Dmix_{(B)}(X)^{t=0}$ and $\Proje \D_{(B)}(X)^{t=0}$ can be described in terms of Soergel modules as well.

First, Soergel's \emph{Endomorphismensatz} states that there is an isomorphism of graded algebras
$$C\cong\Hom{\D_{(B)}(X)}{\Pc_{w_0}}{\Pc_{w_0}}=\bigoplus_{n\in \Z}\Hom{\Dmix_{(B)}(X)}{\gPc_{w_0}}{\gPc_{w_0}(n)}.$$
In \cite{Soe90} this statement is originally proven representation-theoretically for category $\pazocal O.$ There is also a topological proof, due to Bezrukavnikov--Riche, see \cite{BR}. 

Then, Soergel's \emph{Struktursatz} shows that the functors
\begin{align*}
\widehat{\mathbb V}: \Dmix_{(B)}(X)\to C\modules^\Z,\, &M\mapsto\bigoplus_{n\in\Z}\Hom{\Dmix_{(B)}(X)}{\gPc_{w_0}}{M(n)}\\
\mathbb V: D(X)\to C\modules,\, &M\mapsto\Hom{\D_{(B)}(X)}{\Pc_{w_0}}{M}
\end{align*}
are fully faithful on projective perverse objects. In fact there are isomorphisms
\begin{align*}
\widehat{\mathbb V}(\gPc_w)&\cong \mathbb V(\Pc_w)\cong D_{w}.
\end{align*}
Hence there are equivalences of categories
\begin{align*}
\Proje \Dmix_{(B)}(X)^{t=0}&=C\Smodules^\Z\text{ and}\\
\Proje D(X)^{t=0}&=C\Smodules.
\end{align*}
\subsection{Koszul duality} The existence of the following Koszul duality functor $\gKos$ for $\Dmix_{(B)}(X)$ was first conjectured by Beilinson--Ginzburg in \cite{BG} and proven by Soergel in \cite{Soe90} using the combinatorial descriptions in terms of Soergel modules from above. The very elegant formulation using motivic sheaves is due to Soergel--Wendt, \cite{SoeWe}. The functor $\gKos$ can be constructed as the composition
\begin{align*}\gKos: \Dmix_{(B)}(X)&\cong \Hot^b((\Proje \Dmix_{(B)}(X)^{t=0}) \\
&\cong\Hot^b(C\Smodules^\Z)\\
&\cong\Hot^b(\Dmix_{(B^\vee)}(X^\vee)^{w=0})\cong \Dmix_{(B^\vee)}(X^\vee)
\end{align*}
Under this equivalence the projective perverse motivic sheaf $\gPc_{w}$ is sent to the intersection complex $\gEc_w.$ It also intertwines the grading shifts $(n)$ and $(n)[2n].$ For further properties we refer to \cite{BGS}.

We can now consider the \emph{ungraded version} of Koszul duality in exactly the same way, namely, we have equivalences
\begin{align*}\Kos: \D_{(B)}(X)&\cong \Hot^b((\Proje \D_{(B)}(X)^{t=0}) \\
&\cong\Hot^b(C\Smodules)\\
&\cong\Hot^b(\KM_{(B^\vee)}(X^\vee)^{w=0})\cong \KM_{(B^\vee)}(X^\vee)
\end{align*}
Under this equivalence the projective perverse sheaf $\Pc_{w}$ is sent to the \emph{intersection $K$-theory compex} $\mathcal{IK}_w.$  The functor $\Kos$ inherits all the nice properties of $\gKos.$

Combining everything, we hence obtain the quite satisfying commutative diagram
\begin{center}
\begin{tikzcd}
\Dmix_{(B)}(X) \arrow{r}{\gKos} \arrow{d}{v}
& \Dmix_{(B^\vee)}(X^\vee) \arrow{d}{\iota} \\
\D_{(B)}(X)  \arrow{r}{\Kos}
& \KM_{(B^\vee)}(X^\vee).
\end{tikzcd}\end{center}

%% file: Inhalt/tandweight.tex
\label{sec:tandweightstructures}
For the convenience of the reader, we briefly recall the definitions and gluing of $t$-structures and weight structures of triangulated categories.
\subsection{Definitions}
\begin{definition}\cite[Definition 1.3.1]{BBD}
	Let $\Cc$ be a triangulated category. A \emph{$t$-structure} $t$ on $\Cc$ is a pair $t=(\Cc^{t\leq 0},\Cc^{t\geq 0})$ of full subcategories of $\Cc$ such that with $\Cc^{t\leq n}:=\Cc^{t\leq 0}[-n]$ and $\Cc^{t\geq n}:=\Cc^{t\geq 0}[-n]$ the following conditions are satisfied:
	\begin{enumerate}
		\item $\Cc^{t\leq 0}\subseteq \Cc^{t\leq 1}$ and $\Cc^{t\geq 1}\subseteq \Cc^{t\geq 0};$
		\item for all $X\in \Cc^{t\leq 0}$ and $Y\in\Cc^{t\geq 1}$, we have $\Hom{\Cc}{X}{Y}=0;$
		\item for any $X\in \Cc$ there is a distinguished triangle 
		\begin{center}\disttriangle{A}{X}{B}\end{center} with $A\in \Cc^{t\leq 0}$ and $B\in \Cc^{\geq 1}.$
	\end{enumerate}
	The full subcategory $\Cc^{t=0}=\Cc^{t\leq 0}\cap\Cc^{t\geq 0}$ is called the \emph{heart of the $t$-struture}. The $t$-structure is called \emph{bounded} if $\Cc=\bigcup_n\Cc^{t\geq n}=\bigcup_n\Cc^{t\geq n}.$ \todo{A5}
\end{definition}
\begin{definition}\cite[Definition 1.1.1]{Bon}
	Let $\Cc$ be a triangulated category. A \emph{weight structure} $w$ on $\Cc$ is a pair $w=(\Cc^{w\leq 0},\Cc^{w\geq 0})$ of full subcategories of $\Cc,$ which are closed under direct summands, such that with $\Cc^{w\leq n}:=\Cc^{w\leq 0}[-n]$ and $\Cc^{w\geq n}:=\Cc^{w\geq 0}[-n]$ the following conditions are satisfied:
	\begin{enumerate}
		\item $\Cc^{w\leq 0}\subseteq \Cc^{w\leq 1}$ and $\Cc^{w\geq 1}\subseteq \Cc^{w\geq 0};$
		\item for all $X\in \Cc^{w\geq 0}$ and $Y\in\Cc^{w\leq -1}$, we have $\Hom{\Cc}{X}{Y}=0;$
		\item for any $X\in \Cc$ there is a distinguished triangle 
		\begin{center}\disttriangle{A}{X}{B}\end{center} with $A\in \Cc^{w\geq 1}$ and $B\in \Cc^{w\leq 0}.$
	\end{enumerate}
	The full subcategory $\Cc^{w=0}=\Cc^{w\leq 0}\cap\Cc^{w\geq 0}$ is called the \emph{heart of the weight struture}. The weight structure is called \emph{bounded} if $\Cc=\bigcup_n\Cc^{w\geq n}=\bigcup_n\Cc^{w\geq n}.$\todo{A5}
\end{definition}
\begin{remark}\label{rem:tstructuresonderivedcategories}(1) The standard example of a bounded $t$-structure is the bounded derived category $\Der^b(\Aa)$ of an abelian category $\Aa,$ where we set
\begin{align*}
\Der^b(\Aa)^{t\leq 0}&=\setbuild{X\in\Der(\Aa)}{\Hh^i{X}=0 \text{ for all } i>0}\text{ and}\\
\Der^b(\Aa)^{t\geq 0}&=\setbuild{X\in\Der(\Aa)}{\Hh^i{X}=0 \text{ for all } i<0}.
\end{align*}
The standard example of a bounded weight structure is the bounded homotopy category of chain complexes $\Hot^b(\Aa)$ of an \emph{additive} category $\Aa,$ where we set\todo{A5}
\begin{align*}
\Hot^b(\Aa)^{w\leq 0}&=\genbuild{X\in\Hot^b(\Aa)}{X^i=0 \text{ for all } i>0}_{\cong},\\
\Hot^b(\Aa)^{w\geq 0}&=\genbuild{X\in\Hot^b(\Aa)}{X^i=0 \text{ for all } i<0}_{\cong}
\end{align*}
and by $\cong$ we denote closure under isomorphism.
This already showcases an important distinction between $t$-structures and weight structures. While the heart of a $t$-structure is abelian, the heart of a weight structure is only additive in general, and behaves more like the subcategory of projectives or injectives in an abelian category.\\\noindent
(2) We use the cohomological convention for weight and $t$-structures. One can easily translate to the homological convention, by setting $\Cc_{w\leq0}=\Cc^{w\geq0}$ and $\Cc_{w\geq0}=\Cc^{w\leq0}.$
\end{remark}
\begin{proposition}\label{prop:tweighstructures} Let $\Cc$ be a triangulated category with a $t$-structure or weight structure.
The categories $\pazocal D=\Cc^{t\leq0}$, $\Cc^{t\geq0}$, $\Cc^{w\leq0}$ and $\Cc^{w\geq0}$ are extension stable. That is, for any distinguished triangle in $\Cc$
\begin{center}
\disttriangle{A}{B}{C}
\end{center}
with $A,C\in \pazocal D,$ also $B\in \pazocal D.$
\end{proposition}
We will use standard terminology for exactness of functors.
\begin{definition}
Let $F:\Cc_1\to \Cc_2$ between two triangulated categories with t-structures (weight structures). We say that $F$ is left $t$-exact (or left $w$-exact) if $F(\Cc_1^{t\leq 0})\subset \Cc_2^{t\leq 0}$ (or $F(\Cc_1^{w\leq 0})\subset \Cc_2^{w\leq 0}$) and $F$ is right $t$-exact (or left $w$-exact) if $F(\Cc_1^{t\geq 0})\subset \Cc_2^{t\geq 0}$ (or $F(\Cc_1^{w\geq 0}\subset \Cc_2^{w\geq 0}$). We say that $F$ is $t$-exact ($w$-exact) if $F$ is both left and right $t$-exact (or $w$-exact).
\end{definition}
\subsection{Gluing}
As explained in \cite{BBD}, $t$-structures can be glued together. In fact the axiomatic setup required to perform such a gluing also works for weight structures. But there is \emph{subtle and essential difference} in the definition of the gluing of $t$-structures and weight structures, exchanging $*$ and $!$ functors.
\begin{definition}\cite[Section 1.4.3]{BBD} We call sequence of triangulated functors and categories
$$\Cc_Z\stackrel{i_*=i_!}{\to}\Cc\stackrel{j^*=j^!}{\to}\Cc_U$$
a \emph{gluing datum} if the following properties are fulfilled.
\begin{enumerate}
\item The functor $i_*=i_!$ admits triangulated left and right adjoints, denoted by $i^*$ and $i^!.$
\item The functor $j^*=j^!$ admits triangulated left and right adjoints, denoted by $j_!$ and $j_*.$
\item One has $j^*i_*=0.$
\item For all $K\in \Cc$ the units and counits of the adjunctions can be completed to distinguished triangles
\begin{center}
\disttriangle{j_!j^!K}{K}{i_*i^*K}\\
\disttriangle{i_!i^!K}{K}{j_*j^*K}
\end{center}
\item The functors $i_*=i_!,$ $j_!$ and $j^*=j^!$ are fully faithful.
\end{enumerate}
\end{definition}
\begin{theorem}\label{thm:glueing} Assume that $\Cc_Z\stackrel{i_*=i_!}{\to}\Cc\stackrel{j^*=j^!}{\to}\Cc_U$ is a gluing datum.
\begin{enumerate}\item If $(\Cc^{t\leq 0}_U,\Cc^{t\geq 0}_U)$ and $(\Cc^{t\leq 0}_Z,\Cc^{t\geq 0}_Z)$ are $t$-structures on $\Cc_U$ and $\Cc_Z,$ then 
\begin{align*}
\Cc^{t\leq 0}&:=\setbuild{X\in \Cc}{j^!K\in \Cc^{t\leq 0}_U\text{ and }i^*K\in \Cc^{t\leq 0}_Z}\text{ and}\\
\Cc^{t\geq 0}&:=\setbuild{X\in \Cc}{j^*K\in \Cc^{t\geq 0}_U\text{ and }i^!K\in \Cc^{t\geq 0}_Z}
\end{align*}
defines a $t$-structure on $\Cc.$
\item If $(\Cc^{w\leq 0}_U,\Cc^{w\geq 0}_U)$ and $(\Cc^{w\leq 0}_Z,\Cc^{w\geq 0}_Z)$ are weight structures on $\Cc_U$ and $\Cc_Z,$ then 
\begin{align*}
\Cc^{w\leq 0}&:=\setbuild{X\in \Cc}{j^*K\in \Cc^{w\leq 0}_U\text{ and }i^!K\in \Cc^{w\leq 0}_Z}\text{ and}\\
\Cc^{w\geq 0}&:=\setbuild{X\in \Cc}{j^!K\in \Cc^{w\geq 0}_U\text{ and }i^*K\in \Cc^{w\geq 0}_Z}
\end{align*}
defines a weight structure on $\Cc.$
\end{enumerate}
\end{theorem}\begin{proof} The statement for $t$-structures is \cite[Theorem 1.4.10]{BBD}. The statement for weight structures is \cite[Theorem 8.2.3]{Bon}.\end{proof}
\subsection{Weight complex and realisation functors}
It is often possible to realise a triangulated category with $t$-structure as the derived category of its heart. Similarly, one can often realise a triangulated category with a weight structure as the homotopy category of chain complexes of its heart. We recall some statements from the literature.
\begin{theorem}\label{thm:weightcomplex} Let $\Cc$ be an ``enhanced'' triangulated category, meaning that either
\begin{enumerate}
\item (Derivator) $\Cc=\mathbb D (\operatorname{pt})$, where $\mathbb D$ is a strong stable derivator.
\item ($\infty$-category) $\Cc=\operatorname{Ho}(\Cc')$, where $\Cc'$ is a stable $\infty$-category.
\item ($f$-category) There is an $f$-category $\pazocal DF$ over $\Cc.$
\end{enumerate}
Assume that $\Cc$ is equipped with a $t$-structure. Then there is a triangulated functor called \emph{realisation functor}
$$\Der^b(\Cc^{t=0})\rightarrow \Cc $$
restricting to the inclusion of the heart $\Cc^{t=0}\rightarrow \Cc.$

Assume that $\Cc$ is equipped with a bounded weight structure. Then there is a triangulated functor called \emph{weight complex functor}\todo{A5}
$$\Cc\rightarrow\Hot^b(\Cc^{w=0}) $$
restricting to the inclusion of the heart $\Cc^{w=0}\rightarrow \Hot^b(\Cc^{w=0}).$
\end{theorem}
\begin{proof} For the statement about $t$-structures, we refer to \cite{Vir} for derivators, \cite{Lur} for $\infty$-categories and \cite{Bei} for $f$-categories.
For the statement about weight structures, we refer to \cite{Bon} for $f$-categories and \cite{Sos}, \cite{Ko} for $\infty$-categories. In fact, the derivator assumption implies the $f$-category assumption by \cite{Mod}.
\end{proof}
There are different assumptions under which the above functors can be shown to be fully faithful. We refer to the references in the proof above. Furthemore, it can be shown that realisation and weight complex functors are compatible with ``enhanced'' exact triangulated functors between ``enhanced'' triangulated categories. We note that the categories of motives and the six operations between them are all ``enhanced''.